\documentclass[12pt]{article}
\usepackage{amssymb}
\usepackage{amsmath}

\newtheorem{theo}{Theorem}

\newtheorem{lemm}{Lemma}
\newtheorem{coro}{Corollary}
\newtheorem{rema}{Remark}
\newtheorem{Defi}{Definition}

\newcommand{\cqfd}
{%
\mbox{}%
\nolinebreak%
\hfill%
\rule{2mm}{2mm}%
\medbreak%
\par%
}
\newfont{\gothic}{eufb10}
\date{\empty}
\begin{document}
\title{A geometric application of Nori's connectivity theorem}
\author{Claire Voisin\\ Institut de math{\'e}matiques de Jussieu, CNRS,UMR
7586} \maketitle \setcounter{section}{-1}
\section{Introduction}
Our purpose in this paper is to study rational maps from varieties
with small dimensional moduli space to general hypersurfaces in
projective space. In the last section, we shall eventually extend
this to the study of correspondences instead of rational maps.
\begin{Defi}\label{defi}Let $\mathcal{Y}\rightarrow \mathcal{S}$ be a family
of $r$-dimensional smooth projective varieties. We say that a
$d$-dimensional variety $X$ is rationally swept out by  varieties
parametrized by $\mathcal{S}$, if there exist a quasiprojective
variety $B$ of dimension $d-r$, a family ${K}\rightarrow B$ which
is the pull-back of the family $\mathcal{Y}$ via a morphism
$\psi:B\rightarrow\mathcal{S}$, and a dominant rational map
$$\phi:{K}\dashrightarrow X,$$
(which is necessarily generically finite on the generic fiber
$K_b$ since $dim\,{K}=dim\,X$).
\end{Defi}
Our main result in this paper concerns the problem of sweeping out
general hypersurfaces of degree $N\geq d+2$ in $\mathbb{P}^{d+1}$
:
\begin{theo} \label{mainintro} Fix an integer $1\leq r\leq d$. Let $\gamma=\frac{r-1}{2},\,r$
odd, or $\gamma= \frac{r}{2},\,r$ even, that is $\gamma$ is the
round-up of $\frac{r-1}{2}$. Let $\mathcal{Y}\rightarrow
\mathcal{S}$, $dim\,\mathcal{S}=C$, be a family
  of
$r$-dimensional smooth projective varieties. Then the general
hypersurface of degree $N$ in $\mathbb{P}^{d+1}$ is not rationally
swept out by varieties parameterized by $\mathcal{S}$ if
\begin{eqnarray} \label{estimintro}(N+1)r\geq2d+C+2,\,(\gamma+1)N\geq 2d-r+1+C
.\end{eqnarray}
\end{theo}
(Note that except for $r=1$, the second inequality implies the
first.)
\begin{rema} One could of course prove a similar statement for
sufficiently ample hypersurfaces in any smooth variety. In the
case of projective space, the estimates on $N$ are sharp, and
allow applications to the Calabi-Yau case (see section \ref{3}).
\end{rema} The proof is Hodge theoretic. Unlike \cite{clemensens},
\cite{voisinjdg}, \cite{ein}, \cite{chian},  the result has
nothing to do with the canonical bundle of the varieties
$Y_t,\,t\in\mathcal{S}$. Instead, the key point is the fact that
the dimension of the moduli space $\mathcal{S}$ is small :  if
every $X$ was rationally swept out by varieties parameterized by
$\mathcal{S}$, for fixed $Y$, there would be a generically finite
rational map
$$Y\times \widetilde{U}_Y\dashrightarrow
\mathcal{X}_{\widetilde{U}_Y}$$ where
$\mathcal{X}_{\widetilde{U}_Y} $ is the pull-back via a morphism
$\rho: \widetilde{U}_Y\rightarrow U$ of the universal hypersurface
parameterized by $U\subset
H^0(\mathcal{O}_{\mathbb{P}^{d+1}}(N))$, and $Im\,\rho$  is of
codimension $\leq dim\,\mathcal{S}$. This will be shown to
contradict Nori's connectivity theorem (Theorem \ref{connect}).

 We shall apply this
particularly to the case of Calabi-Yau hypersurfaces. In
 the paper
\cite{lang}, Lang formulates a number of conjectures concerning
smooth projective complex varieties $X$. One of them is that the
analytic closure of the union of the images of holomorphic maps
from $\mathbb{C}$ to $X$ is equal to the union of the images of
non constant rational maps from an abelian variety to $X$. Another
one is that this locus is equal to $X$ itself if and only if $X$
is not of general type.

Next, by a standard countability argument for Chow varieties, we
see that, according to these conjectures, if $X$ is not of general
type, there should exist a quasiprojective variety $B$, a family
${K}\rightarrow B$ of abelian varieties, and a dominant rational
map
$$\phi:{K}\dashrightarrow X,$$
which is non constant on the generic fiber $K_b,\,b\in B$.

Let us now consider the case where $X$ is a Calabi-Yau variety,
that is $K_X$ is trivial. We claim that if a map $\phi$ as above
exists, then we may assume that $\phi_{\mid K_b}$ is generically
finite, for generic $b\in B$. Indeed, because $H^0(X,K_X)\not=0$,
for generic $b\in B$, the image $\phi(K_b)$ has effective
canonical bundle, in the sense that any desingularization of it
has effective canonical bundle, as follows from adjunction formula
and the fact that the $\phi(K_b)$ cover $X$. Now it is immediate
to prove that any dominant rational map
$$K_b\dashrightarrow K'_b,$$
where $K_b$ is an abelian variety and $K'_b$ has effective
canonical bundle, factors through the quotient map
 $K_b\rightarrow K''_b$, where $K''_b$ is an abelian
variety, which is a quotient of $K_b$, and has the same dimension
as $K'_b$. Replacing the family of abelian varieties
$({K_b})_{b\in B}$, by the family $(K''_b)_{b\in B}$ gives the
desired $\phi'$.

In other words, Lang's conjecture asserts in particular that a
Calabi-Yau variety should be rationally swept out by
$r$-dimensional abelian varieties, for some $r\geq1$. Our theorem
\ref{mainintro} implies :
\begin{theo}\label{main} Let $X$ be a general Calabi-Yau
hypersurface in projective space $ \mathbb{P}^{d+1}$, that is
$N=d+2$. Then $X$ is not rationally swept out by $r$-dimensional
abelian varieties, for any $r\geq2$.
\end{theo}
Hence, if Lang's conjecture is true, such an $X$ should be swept
out by elliptic curves.

On the other hand, we also prove the following
\begin{lemm} \label{lemmintro}If $X$ is a general Calabi-Yau hypersurface of dimension
$\geq2$, $X$ is not rationally swept out by elliptic curves of
fixed modulus.
\end{lemm}
By ``rationally swept out by elliptic curves of fixed modulus'',
we mean that the elliptic curves in the family ${K}\rightarrow B$
of definition \ref{defi} have constant modulus.

Hence, combining theorem \ref{main} with the above lemma, we get
the following corollary, which was pointed out to us by J. Harris
:
\begin{coro}\label{corintro} If Lang's conjecture is true, any Calabi-Yau
hypersurface $X$ of dimension $\geq2$ has a divisor which is
uniruled.
\end{coro}
 In dimension $3$, this shows that Lang's conjecture and
Clemens conjecture on the finiteness of rational curves of fixed
degree in a general quintic threefold, contradict.

In the case of hypersurfaces of general type, inequality
(\ref{estimintro}) can be applied to give a non trivial estimate
on the minimal genus of covering families of curves, but the
estimate is not sharp and could be obtained directly by geometry.
What is interesting however is that looking more precisely at the
proof of Theorem \ref{mainintro}, we shall see that the result
concerns in fact only the Hodge structure on $H^d(X)_{prim}$ and
not the effective geometry of $X$. In fact we get as well :
\begin{theo}\label{theointrobis} Let $X$ be a general hypersurface of
degree $N\geq 2d-2+3g,\,g\geq2$  or $N\geq 2d+2,\,g=1$, in
$\mathbb{P}^{d+1}$. Then there exists no non-zero morphism of
Hodge structure
$$H^d(X,\mathbb{Q})_{prim}\rightarrow H^d(Y,\mathbb{Q}),$$ where $Y$ is rationally
swept out by curves of genus $g$.
\end{theo}
 Combining this statement
with the generalization of Mumfords theorem \cite{mumford}, this
implies in particular that for $N\geq 2d+2$, $X$ general, there
exists no correspondence $\Gamma\in CH^d(Y\times X)$ inducing a
surjective map
$$\Gamma_*:CH_0(Y)_0 \rightarrow CH_0(X)_0,$$ where $Y$ admits an
elliptic fibration. Similarly, if $g\geq 2$ and $N\geq 2d+2$,
there exists no such correspondence $\Gamma\in CH^d(Y\times X)$
 where $Y$ admits a fibration whose generic fiber is a genus $g$ curve.
  One may wonder whether these statement are true for any
such hypersurface or only for the general one.

 The paper is organized as follows : in
section \ref{1}, we recall briefly the proof of Nori's
connectivity theorem for hypersurfaces in projective space, in
order to extend it to families of hypersurfaces parameterized by
subvarieties of the moduli space which are of small codimension.
This will show us that for any family of  hypersurfaces
parameterized by a subvariety of the moduli space which is of
small codimension, the Hodge level of the cohomology groups of the
total space of the family is small.

 The next
section is devoted to the proof (by contradiction) of theorem
\ref{mainintro}.

 In section \ref{3}, we prove  the applications of this result described above.

\vspace{0,5cm}

{\bf Acknowledgements.} I would like to thank J. Harris who
started me thinking to these problems and to Herb Clemens for very
interesting exchanges.

This work has been essentially done at the University La Sapienza,
and I would like to thanks the organizers of the trimester
``Moduli spaces, Lie Theory, interactions with Physics'', for the
excellent working conditions and atmosphere I found  there.
\section{\label{1}  Nori's connectivity theorem for
hypersurfaces} In this section, we summarize the main points of
the proof of  Nori's connectivity theorem for hypersurfaces, in
order to prove theorem \ref{connect}, which is the precision of it
that we will need. in \cite{AS}, \cite{O}, a sharper study of
similar explicit bounds can be found.

We consider hypersurfaces of degree $N$ in $\mathbb{P}^{d+1}$, and
we assume that $N\geq d+2$. Fix an integer $r$ such that $1\leq
r\leq d$, and let $\gamma$ be the round-up of $\frac{r-1}{2}$.
Denote by $U\subset H^0(\mathcal{O}_{\mathbb{P}^{d+1}}(N))$ the
open set parametrizing smooth hypersurfaces. Let $\rho:
\mathcal{M}\rightarrow U$ be a morphism, where $\mathcal{M}$ is
smooth quasi-projective. We assume that $Corank\,\rho$ is constant
equal to $C$.  We also assume for simplicity that $Im\,\rho$ is
stable under the action of $Gl(N)$. Let $\mathcal{X}_U$ be the
universal hypersurface parametrized by $U$ and
$$\mathcal{X}_{\mathcal{M}}:=\mathcal{X}_U\times_U\mathcal{M}.$$
Let
$$j:\mathcal{X}_{\mathcal{M}}\hookrightarrow
\mathcal{M}\times{\mathbb{P}^{d+1}}$$ be the natural embedding.
$\mathcal{X}_{\mathcal{M}}$ is a smooth quasi-projective variety,
hence its cohomology groups carry mixed Hodge structures with
associated Hodge filtration
$F^iH^k(\mathcal{X}_{\mathcal{M}},\mathbb{C})$.
\begin{theo}\label{connect} i) Assume that
\begin{eqnarray}\label{estimate}
(N+1)r \geq 2d+C+2  . \end{eqnarray} Then,  the restriction map
$$j^*:F^{d}H^{2d-r}(\mathcal{M}\times{\mathbb{P}^{d+1}},\mathbb{C})
\rightarrow F^{d}H^{2d-r}(\mathcal{X}_{\mathcal{M}},\mathbb{C})$$
is surjective.

ii) If
\begin{eqnarray}\label{estimate1}
(\gamma+1)N\geq 2d+1-r+C,
\end{eqnarray} then for any $i\geq1$, the restriction map
$$j^*:H^{2d-r-i}(\mathcal{M}\times{\mathbb{P}^{d+1}},\mathbb{C})
\rightarrow H^{2d-r-i}(\mathcal{X}_{\mathcal{M}},\mathbb{C})$$ is
surjective.
\end{theo}  {\bf Proof.} i)
One first reduces   i), (see \cite{nori}), to proving that under
the assumption (\ref{estimate}), the restriction map
$$j^*:H^l(\Omega_{\mathcal{M}\times{\mathbb{P}^{d+1}}}^k)\rightarrow
H^l(\Omega_{\mathcal{X}_{\mathcal{M}}}^k)$$ is bijective, for
$l\leq d-r,\,k+l\leq 2d-r$. This step uses the mixed Hodge
structure on relative cohomology and the Fr\"ohlicher spectral
sequence.

Denote respectively by $\pi_X,\,\pi_{\mathbb{P}}$ the natural maps
$$\mathcal{X}_{\mathcal{M}}\rightarrow \mathcal{M},\,
\mathcal{M}\times{\mathbb{P}^{d+1}}\rightarrow \mathcal{M}.$$ A
Leray spectral sequence argument shows that it suffices to prove
that under the assumption (\ref{estimate}) one has :

\begin{eqnarray} \label{numero}  The\,\, restriction\,\,
map\,\,  j^*:R^l\pi_{{\mathbb
P}*}(\Omega_{\mathcal{M}\times{\mathbb{P}^{d+1}}}^k)\rightarrow
R^l\pi_{X*}(\Omega_{\mathcal{X}_{\mathcal{M}}}^k) \,\, is\,\,
bijective
\end{eqnarray}
 $for \,\,l\leq d-r,\,k+l\leq 2d-r$.

Let
$$\mathcal{H}^d_{prim},\,\mathcal{H}^{p,q}_{prim},\,p+q=d,$$ be the Hodge
bundles associated to the variation of Hodge structure on the
primitive cohomology of the family
$\pi_X:\mathcal{X}_{\mathcal{M}}\rightarrow \mathcal{M}$. The
infinitesimal variation of Hodge structure on the primitive
cohomology of the fibers of $\pi_X$ is described by maps
$$\overline{\nabla}:\mathcal{H}^{p,q}_{prim}
\rightarrow\mathcal{H}^{p-1,q+1}_{prim}\otimes\Omega_{\mathcal{M}},$$
and they can be iterated to produce a complex :
\begin{eqnarray}\label{drdol}
\ldots
\mathcal{H}^{p+1,q-1}_{prim}\otimes\Omega_{\mathcal{M}}^{s-1}
\stackrel{\overline{\nabla}}{\rightarrow}
\mathcal{H}^{p,q}_{prim}\otimes\Omega_{\mathcal{M}}^{s}
\stackrel{\overline{\nabla}}{\rightarrow}
\mathcal{H}^{p-1,q+1}_{prim}\otimes\Omega_{\mathcal{M}}^{s+1}\ldots
\end{eqnarray}
 One can show, using the filtration of $\Omega_{\mathcal{X}_{\mathcal{M}}}^k$
 by the subbundles
 $\pi_X^*\Omega_B^s\wedge\Omega_{\mathcal{X}_{\mathcal{M}}}^{k-s}$and
 the associated spectral sequence,
  that (\ref{numero}) is equivalent to the following

\vspace{0,5cm}

 {\it The sequence (\ref{drdol}) is exact at the
middle
 for $q\leq d-r,\,p+s+q\leq2d-r$.}

\vspace{0,5cm}

 Note that since $p+q=d$, the last inequality reduces to
$s\leq d-r$.

It is convenient to dualize (\ref{drdol}) using Serre duality,
which gives :
\begin{eqnarray}
{\label{dual}}
\mathcal{H}^{q+1,p-1}_{prim}\otimes\bigwedge^{s+1}T_{\mathcal{M}}
\stackrel{^t\overline{\nabla}}{\rightarrow}
\mathcal{H}^{q,p}_{prim}\otimes\bigwedge^sT_{\mathcal{M}}^{s}
\stackrel{^t\overline{\nabla}}{\rightarrow}
\mathcal{H}^{q-1,p+1}_{prim}\otimes\bigwedge^{s-1}T_{\mathcal{M}}^{s+1}.
\end{eqnarray}
We finally use Griffiths, Griffiths-Carlson description of the
IVHS of hypersurfaces (\cite{griffiths}, \cite{carlson}) to
describe the complex (\ref{dual}) at the point $f\in \mathcal{M}$
as follows. We have the map $\rho_*:T_{\mathcal{M},f}\rightarrow
T_{U,f}= S^{N}$, where $S$ is the polynomial ring in $d+2$
variables. Next the residue map provides  isomorphisms
$$  R_f^{-d-2+N(p+1)}\cong H^{q,p}_{prim}(X_f),$$
where $R_f:=S/J_f$ is the Jacobian ideal of $f$, and $R_f^k$
denotes its degree $k$ component. The map $\overline{\nabla}$
identifies then, up to a coefficient, to the map given by
multiplication
$$
R_f^{-d-2+N(p+1)}\rightarrow
Hom\,(T_{\mathcal{M},f},R_f^{-d-2+N(p+2)}).$$ It follows from this
that the sequence (\ref{dual}) identifies to the following piece
of the Koszul complex of the Jacobian ring $R_f$ with respect to
the action of $T_{\mathcal{M},f}$ on it by multiplication:
\begin{eqnarray}\label{koszul}
R_f^{-d-2+Np}\otimes\bigwedge^{s+1}T_{\mathcal{M},f}
\stackrel{\delta}{\rightarrow}
R_f^{-d-2+N(p+1)}\otimes\bigwedge^sT_{\mathcal{M},f}
\stackrel{\delta}{\rightarrow}
R_f^{-d-2+N(p+2)}\otimes\bigwedge^{s-1}T_{\mathcal{M},f}.
\end{eqnarray}
Now, by assumption, if $W$ is the image of $\rho_*$, $W\subset
S^{N}$ is a base-point free linear system, because it contains the
jacobian ideal $J_f^{N}$, and it satisfies $codim\,W= C$.

One  verifies that it suffices to check exactness at the middle of
the exact sequences (\ref{koszul}) in the considered range, with
$T_{\mathcal{M},f}$ replaced with $W$. This last fact is then a
consequence of the following theorem due to M. Green :
\begin{theo}\label{theogreen}\cite{green} Let $W\subset S^{N}$ be a base-point free
linear system. Then the following sequence, where the
differentials are the Koszul differentials
\begin{eqnarray}\label{koszgre}S^{-d-2+Np}\otimes\bigwedge^{s+1}W
\stackrel{\delta}{\rightarrow} S^{-d-2+N(p+1)}\otimes\bigwedge^sW
\stackrel{\delta}{\rightarrow}
S^{-d-2+N(p+2)}\otimes\bigwedge^{s-1}W
\end{eqnarray}
is exact for $-d-2+Np\geq s+codim\,W$.
\end{theo}
Using the fact that the Jacobian ideal is generated by a regular
sequence in degree $N-1$, one then shows that the same is true
when $S^i$ is replaced with $R_f^i$ in (\ref{koszgre}), at least
if $-d-2+N(p+1)\geq N-1$.

We now conclude the proof of i). We have just proved that
(\ref{drdol}) is exact at the middle if
$$-d-2+Np\geq s+C,\,-d-2+N(p+1)\geq
N-1 .
$$
Since we assumed $N\geq d+2$, the second inequality is satisfied
when $p\geq1$. Next, if $ q\leq d-r,\,s\leq d-r$, we have
$$p\geq r\geq1,\,s\leq d-r.$$
Hence, the exactness of (\ref{drdol}) in the range $ q\leq
d-r,\,s\leq d-r$ will follow from the inequality
$$-d-2+Nr\geq d-r+C,$$
that is (\ref{estimate}).

\vspace{0,5cm}

ii) The proof is exactly similar, and we just sketch it in order
to see where the numerical assumption is used. We first observe
that by a mixed Hodge structure argument (cf \cite{nori}), it
suffices, in order to get the surjectivity of the restriction map
:
$$j^*:H^{2d-r-i}(\mathcal{M}\times{\mathbb{P}^{d+1}},\mathbb{C})
\rightarrow H^{2d-r-i}(\mathcal{X}_{\mathcal{M}},\mathbb{C}),$$ to
show the surjectivity of the restriction map :
$$j^*:F^{d-r+\gamma_i}H^{2d-r-i}(\mathcal{M}\times{\mathbb{P}^{d+1}},\mathbb{C})
\rightarrow
F^{d-r+\gamma_i}H^{2d-r-i}(\mathcal{X}_{\mathcal{M}},\mathbb{C}),$$
where $\gamma_i$ is the round-up of $\frac{r-i}{2}$. (This is
because the round-up of $\frac{2d-r-i}{2}$ is $d-r+\gamma_i$.) We
reduce then this last fact to showing :
\begin{eqnarray} \label{numero1}  The\,\, restriction\,\,
map\,\,  j^*:R^l\pi_{{\mathbb
P}*}(\Omega_{\mathcal{M}\times{\mathbb{P}^{d+1}}}^k)\rightarrow
R^l\pi_{X*}(\Omega_{\mathcal{X}_{\mathcal{M}}}^k) \,\, is\,\,
bijective
\end{eqnarray}
 $for \,\,l\leq d-i-\gamma_i,\,k+l\leq 2d-r-i$.

 Expressing the cohomology groups above with the help of the IVHS
 on the primitive cohomology of the fibers of $\pi_X$, this is
 reduced to proving :

\vspace{0,5cm}

 {\it The sequence (\ref{drdol}) is exact at the
middle
 for $q\leq d-i-\gamma_i,\,p+q=d,\,p+s+q
\leq2d-r-i$.}

\vspace{0,5cm}

 Using the Carlson-Griffiths theory, we are now reduced to prove :

\vspace{0,5cm}

 {\it the following sequence :
\begin{eqnarray}\label{koszul1}
R_f^{-d-2+Np}\otimes\bigwedge^{s+1}T_{\mathcal{M},f}
\stackrel{\delta}{\rightarrow}
R_f^{-d-2+N(p+1)}\otimes\bigwedge^sT_{\mathcal{M},f}
\stackrel{\delta}{\rightarrow}
R_f^{-d-2+N(p+2)}\otimes\bigwedge^{s-1}T_{\mathcal{M},f}.
\end{eqnarray}
is exact for $p\geq \gamma_i+i$, $s\leq d-r-i$.}

\vspace{0,5cm}

 As in the previous proof, we now apply the theorem \ref{theogreen} and
conclude that the last statement is true if \begin{eqnarray}
\label{estimate2}-d-2+N(\gamma_i+i)\geq C+d-r-i. \end{eqnarray}
Now it is clear that the $\gamma_i+i$ are increasing with $i$,
while the $C+d-r-i$ are decreasing with $i$. Hence it suffices to
have (\ref{estimate2}) satisfied for $i=1$, which is exactly
inequality (\ref{estimate1}).
 \cqfd

 Denoting by
 $H^{2d-r}(\mathcal{X}_{\mathcal{M}})_{prim}$ the quotient
 $$H^{2d-r}(\mathcal{X}_{\mathcal{M}})_{prim}/j^*(H^{2d-r}({\mathcal{M}}
 \times{\mathbb{P}^{d+1}})),$$
 we shall only be interested with the pure part
 $$W_{2d-r}H^{2d-r}(\mathcal{X}_{\mathcal{M}})_{prim},$$
 which is the part of the cohomology which comes from any smooth
 projective
compactification of $\mathcal{X}_{\mathcal{M}}$. It carries a pure
Hodge structure of weight $2d-r$.
\begin{coro}\label{coro24mai} Under the assumptions of theorem \ref{connect}, the Hodge structure on
$W_{2d-r}H^{2d-r}(\mathcal{X}_{\mathcal{M}})_{prim}$ is of Hodge
level $\leq r-2$.
\end{coro}
{\bf Proof.} Recall that the Hodge level of a Hodge structure
$H,\,H_{\mathbb{C}}=\oplus H^{p,q}$ is
$$Max\,\{p-q,\,H^{p,q}\not=0\}.$$
Since we know that
$F^{d}W_{2d-r}H^{2d-r}(\mathcal{X}_{\mathcal{M}})_{prim}=0$, we
have
$$H^{p,q}(W_{2d-r}H^{2d-r}(\mathcal{X}_{\mathcal{M}})_{prim})=0,\,
\,for\,\, p\geq d.$$ Since $H^{p,q}=0$ for $p+q\not=2d-r$, it
follows that the Hodge level is $\leq d-1-(2d-r-(d-1))=r-2$. \cqfd

\section{\label{2}Proof of theorem \ref{mainintro}}
We prove theorem \ref{mainintro} by contradiction. Using Chow
varieties, or relative Hilbert schemes, we see that there exist
countably many quasi-projective varieties $\mathcal{B}$
parameterizing triples $(t,f,\phi_s)$, where $t\in\mathcal{S}$,
$f\in U$, and $\phi_s$ is a rational map
$\phi_s:Y_s\dashrightarrow X_f$ which is generically finite onto
its image. For fixed generic $f$, our assumption is that the
images of such $\phi_s$ fill-in $X_f$, and a countability argument
then shows  that there exists one $\mathcal{B}$, which dominates
$U$ via the second projection, and which is such that the
universal rational map
$$\Phi:\mathcal{K}\dashrightarrow \mathcal{X}_U$$
is dominating, where
$$\pi:\mathcal{K}\rightarrow\mathcal{B}$$
is the pull-back via the first projection
$\Psi:\mathcal{B}\rightarrow\mathcal{S}$ of the family
$\mathcal{Y}\rightarrow\mathcal{S}$, and, as in the previous
section, $\mathcal{X}_U$ is the universal hypersurface
parameterized by $U$. We shall denote by $B_f$ the (generic) fiber
of the second projection $q:\mathcal{B}\rightarrow U$ and
$\pi_f:K_f\rightarrow B_f$ the induced family.  By taking
desingularizations, we may assume that $ \mathcal{B}$ hence
$\mathcal{K}$ are smooth, and by assumption the map $\pi$ is
smooth. Since $f$ is generic, $B_f$ and $K_f$ are then also
smooth. Finally, we may, up to replacing $\mathcal{B}$ by a closed
subvariety, assume that the restriction $\phi_f:K_f\dashrightarrow
X_f$ of $\Phi$ to $K_f$ is generically finite and dominating. In
particular, $dim\,B_f=d-r$.

Now we make the following construction : denote by $\mathcal{B}_f$
the space
$$\mathcal{B}\times_{\mathcal{S}}B_f.$$
Restricting to  Zariski open sets of $\mathcal{B}$ and $B_f$, we
may assume that $\mathcal{B}_f$ is smooth. The generic point of
$\mathcal{B}_f$ parameterizes, via the second projection
$p:\mathcal{B}_f\rightarrow B_f$, a variety $Y_t$ together with a
rational map $\phi_{t,f}:Y_t\dashrightarrow X_f$, and, via the
second projection, a rational map
$$\phi_{t,f'}:Y_t\dashrightarrow X_{f'}.$$
Let $\rho:\mathcal{B}_f\rightarrow U$ be the composition of the
first projection and the map $q:\mathcal{B}\rightarrow U$ and let
$m:\mathcal{B}_f\rightarrow\mathcal{S}$ be the natural map. We
shall also use the notation
$$\mathcal{K}_f=\mathcal{K}\times_{\mathcal{B}}{\mathcal{B}_f}=\bigcup_{t\in \mathcal{B}_f}Y_{m(t)},$$
and $$ \mathcal{X}_f:=
\mathcal{X}\times_U{\mathcal{B}_f}=\bigcup_{t\in
\mathcal{B}_f}X_{\rho(t)}.$$ The map $\Phi$ induces a rational map
$\Phi_f:\mathcal{K}_f\dashrightarrow \mathcal{X}_f$ which is
compatible with the maps $\mathcal{K}_f\rightarrow
{\mathcal{B}_f}$ and $\mathcal{X}_f\rightarrow {\mathcal{B}_f}$.
It follows that the graph $\Gamma$ of $\Phi$ is contained in
$$\mathcal{Y}_f:= \mathcal{K}_f\times_{{\mathcal{B}_f}}\mathcal{X}_f
\cong K_f\times_{B_f}\mathcal{X}_f$$ and is of codimension $d$ in
$\mathcal{Y}_f$.

Note that $\mathcal{Y}_f$ contains $K_f\times X_f$ and that
$\Gamma\cap {K_f\times X_f}$ is nothing but the graph of $\phi_f$.
Now, the class of this last graph in $H^{2d}(K_f\times
X_f,\mathbb{Q})$ does not vanish in
$$H^{2d}(K_f\times X_f,\mathbb{Q})/H^{2d}(K_f\times \mathbb{P}^{d+1},\mathbb{Q}).$$
Indeed, its K\"unneth  component $\phi_f^*$ in
$$Hom\,(H^d(X,\mathbb{Q})_{prim},H^d(K_f,\mathbb{Q}))$$
does not vanish, because $N\geq d+2$, so that the transcendent
part of $H^d(X_f,\mathbb{Q})_{prim}$, that is the orthogonal of
all sub-Hodge structures which are of level $<d$, is non-zero, so
that it cannot be  annihilated by $\phi_f^*$, because $\phi_f$ is
dominating.

It follows that the class $\gamma$ of $\Gamma$ in $H^{2d}(
\mathcal{Y}_f,\mathbb{Q}) $ is non zero modulo
$H^{2d}(\mathcal{K}_f\times \mathbb{P}^{d+1},\mathbb{Q})$.

Recall next that with the help of a polarization, that is a choice
of a relatively ample line bundle on
$\mathcal{Y}_f\stackrel{\pi}{\rightarrow} \mathcal{X}_f$, the
cohomology $H^{2d}( \mathcal{Y}_f,\mathbb{Q}) $ splits into a
direct sum
$$H^{2d}( \mathcal{Y}_f,\mathbb{Q})=\bigoplus_l
\,H^{2d-l}(\mathcal{X}_f,R^l\pi_*\mathbb{Q}).$$ It is easy to
check by similar reasons as above  that the component $\gamma_r$
of $\gamma$ in
$$H^{2d-r}(\mathcal{X}_f,R^r\pi_*\mathbb{Q})/
H^{2d-r}(\mathcal{B}_f\times
\mathbb{P}^{d+1},R^r\pi'_*\mathbb{Q}),$$ where in the second term
$\pi'$ is the natural map
$$\mathcal{K}_f\times\mathbb{P}^{d+1}\rightarrow \mathcal{B}_f\times
\mathbb{P}^{d+1},$$ is non zero. We may finally for the same
reason refine this, replacing $R^r\pi_*\mathbb{Q}$ by its quotient
$R^r\pi_*\mathbb{Q}_{tr}$, that is its quotient by the maximal
sub-Hodge structure which exists generically on $B_f$ and is not
of maximal Hodge level $r$.
 Note that
$\gamma$, $\gamma_r$ and $\gamma_{r,tr}$ are Hodge classes, that
is belong to the $F^dW_{2d}$-level of the considered cohomology
groups, which all have mixed Hodge structures.

Next we denote by
$$g:\mathcal{B}_f\rightarrow B_f,\,\tilde{g}:=g\circ\pi_X:\mathcal{X}_f\rightarrow B_f,$$
the natural maps. We observe that shrinking $B_f$ if necessary,
the fibers $\mathcal{B}_{f,t}$ of $g$ are smooth and the map
$\rho_{\mid \mathcal{B}_{f,t}}$ has constant corank  $\geq
C=dim\,\mathcal{S}$, because $\rho:\mathcal{B}\rightarrow U$ is
submersive near $B_f$, and the fibers of $g$ identify to the
fibers of $\mathcal{B}\rightarrow\mathcal{S}$.

Hence we may apply theorem \ref{connect} and its corollary. It
says that under the assumtions \ref{estimintro}, the fibers
$\mathcal{X}_{f,t}$ of $\tilde{g}$ satisfy
$H^{2d-r-i}(\mathcal{X}_{f,t},\mathbb{Q})=0,\,i\geq1$.

On the other hand, we  observe  that because
$\pi:\mathcal{Y}_f\rightarrow\mathcal{X}_f$ is a fiber product
$K_f\times_{B_f}\mathcal{X}_f$, the local system
$R^r\pi_*\mathbb{Q} $ is a pull-back from $B_f$ :
$$ R^r\pi_*\mathbb{Q}\cong \tilde{g}^{-1}(R^r\pi_{f*}\mathbb{Q}).$$
The vanishing above implies then that
$$(R^{2d-r-i}g_*(R^r\pi_{*}\mathbb{Q}_{tr}))_{prim}=0,\,i\geq1,$$
where $``prim''$ here denotes the quotient by the ambiant
cohomology (while before it was used to mean the orthogonal to
ambiant cohomology). It follows by Leray spectral sequence that
the class $\gamma_{r,tr}$ does not vanish along the fibers of $g$,
that is in $H^0( R^{2d-r}g_*(R^r\pi_{*}\mathbb{Q}_{tr})$ and by
restriction to the general fiber, it does not vanish in
$$H^r(Y_t,\mathbb{Q})_{tr}\otimes H^{2d-r}(\mathcal{X}_{f,t},\mathbb{Q})/
H^{2d-r}(\mathcal{B}_{f,t}\times \mathbb{P}^{d+1},\mathbb{Q}).$$
Note that it is a Hodge class in
$$H^r(Y_t,\mathbb{Q})_{tr})\otimes
H^{2d-r}(\mathcal{X}_{f,t},\mathbb{Q})/H^{2d-r}(\mathcal{B}_{f,t}\times
\mathbb{P}^{d+1},\mathbb{Q}))$$ $$\cong
=Hom\,(H^r(Y_t,\mathbb{Q})_{tr},H^{2d-r}(\mathcal{X}_{f,t},\mathbb{Q})/H^{2d-r}(\mathcal{B}_{f,t}\times
\mathbb{P}^{d+1},\mathbb{Q})).$$  Since the mixed Hodge structure
on the left is pure, this implies that it belongs to
$$Hom_{MHS}(H^r(Y_t,\mathbb{Q})_{tr},H^{2d-r}(\mathcal{X}_{f,t},\mathbb{Q})
/H^{2d-r}(\mathcal{B}_{f,t}\times \mathbb{P}^{d+1},\mathbb{Q}))$$
$$\cong
Hom_{HS}(H^r(Y_t,\mathbb{Q})_{tr},W_{2d-r}H^{2d-r}(\mathcal{X}_{f,t},\mathbb{Q})
/H^{2d-r}(\mathcal{B}_{f,t}\times \mathbb{P}^{d+1},\mathbb{Q}))
,$$ where ``MHS'' and ``HS'' mean   morphisms of mixed
(respectively pure) Hodge structures.

Since  by definition $H^r(Y_t,\mathbb{Q})_{tr}$ has no quotient
Hodge structure which is of Hodge level $<r$, we get now a
contradiction with corollary \ref{coro24mai}, which says that the
Hodge structure on
$W_{2d-r}H^{2d-r}(\mathcal{X}_{f,t},\mathbb{Q})/H^{2d-r}(\mathcal{B}_{f,t}\times
\mathbb{P}^{d+1},\mathbb{Q})$ has Hodge level $<r$.
 \cqfd

\section{\label{3}Rational maps from abelian varieties to Calabi-Yau hypersurfaces
and other applications} Let us apply theorem \ref{mainintro} to
the case of Calabi-Yau hypersurfaces, that is hypersurfaces of
degree $N=d+2$ in $\mathbb{P}^{d+1}$. The moduli space of
$r$-dimensional abelian varieties with given polarization type is
of dimension $\frac{r(r+1)}{2}$. Hence  the conditions of
 theorem \ref{mainintro} become :
 $$(d+3)r\geq2d+\frac{r(r+1)}{2}+2,\,(\gamma+1)(d+2)\geq 2d-r+1+\frac{r(r+1)}{2}
 $$.  It is not hard to check that this is satisfied
 for $2\leq r\leq d$. Hence we get in this case theorem \ref{main}. \cqfd

When $r=1$, the inequality (\ref{estimintro}) is never satisfied
so that
 our argument definitely does not apply to the study of elliptic
 curves in Calabi-Yau hypersurfaces. In fact we could adapt our proof
 of theorem \ref{main} to work as well
 for  Calabi-Yau hypersurfaces in a product of projective spaces. On the other hand,
 certain generic Calabi-Yau  hypersurfaces in a product of projective spaces
 are swept out by elliptic curves,
 eg the hypersurface of bidegree $(3,3)$ in
 $\mathbb{P}^2\times\mathbb{P}^2$. This shows that for
 $r=1$, a different argument has to be found.

  We can however prove the
 following :
\begin{lemm} \label{lemmintrobis}If $X$ is a Calabi-Yau hypersurface of dimension
$\geq2$, $X$ is not rationally swept out by elliptic curves of
fixed modulus.
\end{lemm}
{\bf Proof.} Indeed, fixing otherwise the modulus of the elliptic
curve, we would get, for at least one elliptic curve $E$, an
hypersurface $\mathcal{M}_E$ in the moduli space $\mathcal{M}$ of
$X$ consisting of $X_f$'s which are rationally dominated by some
$E\times B$. For such an $X_f$, there must be an inclusion of
rational Hodge structures induced by the dominant rational map
$\phi:E\times B\dashrightarrow X_f$:
\begin{eqnarray}
\label{morph}\phi^*:H^d(X_f)_{prim}\hookrightarrow H^1(E)\otimes
H^{d-1}(B), \end{eqnarray} because the Hodge structure on
$H^d(X)_{prim}$ is simple.

If we now let  $f$ vary in $\mathcal{M}_E$,  only $B$ deforms with
$f$, not $E$, and it follows that the infinitesimal variation of
Hodge structure on $H^1(E)\otimes H^{d-1}(B)$ $$\overline{\nabla}:
H^{p,q}(H^1(E)\otimes H^{d-1}(B))\rightarrow
H^{p-1,q+1}(H^1(E)\otimes H^{d-1}(B))\otimes
\Omega_{\mathcal{M}_E}$$ has the following form at the point $f\in
\mathcal{M}_E$ : $$ \overline{\nabla}
(\alpha\otimes\beta)=\alpha\otimes\overline{\nabla}_B(\beta)$$ for
$\alpha\in H^{r,s}(E),\,\beta\in H^{p-r,q-s}(B)$, where
$\overline{\nabla}_B$ is the infinitesimal variation of Hodge
structure on $H^{d-1}(B)$. Hence the Yukawa couplings of the IVHS
on $H^1(E)\otimes H^{d-1}(B)$, that is the iterations of
$\overline{\nabla}$, have the following property:

\vspace{0,5cm}

 {\it $\forall\eta\in H^{d,0}(H^1(E)\otimes
H^{d-1}(B))$, the map
$$\overline{\nabla}^d(\eta):S^{d}T_{\mathcal{M}_E,f}\rightarrow
H^{0,d}(H^1(E)\otimes H^{d-1}(B))$$ vanishes.}

\vspace{0,5cm}

 If  there is along ${\mathcal{M}}_E$ an injective morphism of
Hodge structure (\ref{morph}), it follows that the same property
is true for the IVHS of the family of $X_f$'s parameterized by
$\mathcal{M}_E$, namely the Yukawa couplings of $X_f$ vanish on
the hyperplane $K:=T_{\mathcal{M}_E,f}\subset S^{d+2}$. The
Carlson-Griffiths theory \cite{griffiths}, \cite{carlson} shows
easily that this is not the case. Indeed, these Yukawa couplings
identify to the multiplication map
$$S^{d}(S^{d+2})\rightarrow R^{d(d+2)}_f.$$ Assume they vanish on
$K$. Since $K$ is a hyperplane in $S^{d+2}$, the subspace
$$K':=[K:S^1]\subset S^{d+1}$$
has codimension $\leq d+2$. It is without base-point, since
$T_{\mathcal{M}_E,f}$ contains $J_f$. It follows then from
\cite{greengotz}, that $S^{d+3}K'=S^{2d+4}$. But $K^2$ contains
$S^1K'\cdot K =K'\cdot S^1K=W\cdot S^{d+3}=S^{2d+4}$. Hence
$K^2=S^{2d+4}$ and similarly $K^d=S^{d(d+2)}$ contradicting the
fact that $K^d\subset J^{d(d+2)}_f$.
 \cqfd
\begin{rema} In the case of odd dimensional varieties, one can
also use the Mumford-Tate group argument due to Deligne
(\cite{deligne}, p 224) to get this result.
\end{rema}
\begin{coro}\label{corintro1} If Lang's conjecture is true, any Calabi-Yau
hypersurface $X$ of dimension $\geq2$ has a divisor which is
uniruled.
\end{coro}
{\bf Proof.} Indeed, we know that $X$ is rationally swept out by
elliptic curves, but not by elliptic curves with constant modulus.
Hence there is a diagram
$$\begin{matrix}&\widetilde{\overline{\mathcal{K}}}&\stackrel{\phi}{\rightarrow}&X \\
&\pi\downarrow&&\\& \widetilde{\overline{B}}&& , \end{matrix}$$
where we may assume that $\widetilde{\overline{K}}$ and
$\widetilde{\overline{B}}$ are smooth, projective, where
$\widetilde{\overline{\mathcal{K}}}$ is a smooth projective model
of the family $\mathcal{K}\rightarrow B$ on which $\phi$ is
defined, and  that the map
$j:\widetilde{\overline{B}}\rightarrow\mathbb{P}^1$ is defined,
Now, since $\phi$ is generically finite, for generic $t\in
{\mathbb P}^1$, the divisor
$\widetilde{\overline{\mathcal{K}}}_t:=(j\circ\pi)^{-1}(t)$ must
be sent by $\phi$ onto a divisor of $X$, and it follows that for
any $t$ the image by $\phi$ of the divisor
$\widetilde{\overline{\mathcal{K}}}_t:=(j\circ\pi)^{-1}(t)$ must
contain a divisor of $X$. Taking $t=\infty$, and noting that any
component of  $\widetilde{ \overline{\mathcal{K}}}_\infty$ is
uniruled,
 gives
the result. \cqfd

 Finally, we turn to Theorem \ref{theointrobis}. We simply  note
 for this that in the proof of Theorem \ref{mainintro},
we used the dominating rational map
$$\phi:K\dashrightarrow X$$ only to deduce that there is a
corresponding inclusion
$$\phi_f^*:H^d(X_f)_{prim}\rightarrow H^d(K_f)$$ of Hodge structures.
 The Chow argument we used would work
equally with graphs of rational maps replaced with cycles in the
product
$$Y_t\times X.$$
Hence we conclude that everywhere in the paper, we could replace
``dominating rational maps'' by ``cycles in $CH^d(K\times X)$
inducing a non-zero morphism of Hodge structure''
$$H^d(X)_{prim}\rightarrow H^d(K).$$
(Note that such a non-zero morphism should be in fact injective
because the Hodge structure on $H^d(X)_{prim}$ is simple for
generic $X$.) Using the fact that the moduli space of curves of
genus $g\geq1$ has dimension $3g-3,\,g\geq2$ or $1,\,g=1$ , we see
that theorem \ref{theointrobis} is then a consequence of theorem
\ref{mainintro}, where we replace dominating rational maps with
correspondences inducing a non zero morphism of Hodge structure.

\cqfd

\end{document}